\documentclass[11pt]{amsproc}

\usepackage{eufrak}
\usepackage{amssymb}

\theoremstyle{plain}
\newtheorem{thm}{Theorem}[section]
\newtheorem{lem}[thm]{Lemma}
\newtheorem{prop}[thm]{Proposition}

\theoremstyle{remark}
\newtheorem{rem}[thm]{Remark}
\newtheorem{ntn}[thm]{Notation}

\newtheorem*{ta}{\bf Theorem 1}
\newtheorem*{tb}{\bf Theorem 2}

\begin{document}

\title[]{ On a homotopy relation between the
${\mathbf 2}$-local geometry and the Bouc complex for the sporadic group
$\mathbf{Co_3}$}
\maketitle

\begin{center}
\bigskip
{\Large John Maginnis\footnote{Email address: maginnis@math.ksu.edu.} and Silvia Onofrei\footnote{Email address: onofrei@math.ksu.edu.}}\\
{\it Department of Mathematics, Kansas State University,\\ 138 Cardwell
Hall, Manhattan, KS 66506}

\medskip
{\it Journal of Algebra(2007), doi:10.1016/j.algebra.2007.04.028}\\
\end{center}

\bigskip
{\bf Abstract}\\
We study the homotopy relation between the standard $2$-local geometry $\Delta$ and the Bouc complex for the sporadic
group $Co_3$. We also give a result concerning the relative
projectivity of the reduced Lefschetz module $\tilde{L} (\Delta)$.\\


\bigskip
\section{Introduction}
If $G$ is a Lie group in natural characteristic $p$, then the Tits building is the simplicial complex on the maximal parabolic subgroups of $G$. The Bouc complex $|\mathcal{B}_p(G)|$, in this case the complex of unipotent radicals, is the barycentric subdivision of the Tits building.

\medskip
In the $70$'s and $80$'s, Buekenhout \cite{bue}, Ronan and Smith \cite{rs80} and Ronan and Stroth \cite{rst} constructed various geometries for the sporadic simple groups in an attempt to generalize the Tits buildings for Lie groups. Also, Brown \cite{br}, Quillen \cite{qu78}, Bouc \cite{bouc} and others were considering various collections of $p$-subgroups related to group cohomology.

\medskip
In his expository paper, Webb \cite{wbarc} noted the connections between the group geometries and the subgroup complexes, and investigated the associated Lefschetz modules and the corresponding homology decompositions. Webb's paper contains specific calculations, with details from various authors, on geometries and subgroup complexes.

\medskip
The relationship between $p$-local geometries and the Quillen complex of
elementary abelian $p$-subgroups or the
Bouc complex of $p$-radical subgroups for certain sporadic groups was further investigated in papers by Ryba, Smith and Yoshiara \cite{rsy} and Smith and Yoshiara \cite{sy}. In the latter paper, the authors studied the projectivity of the reduced Lefschetz module for some group geometries and observed that the sporadic groups of characteristic $p$-type behave similarly to the finite Lie groups in defining characteristic $p$.

\medskip
There is no uniform approach for $p$-local geometries. A recent, comprehensive study of $2$-local geometries for the $26$ sporadic simple groups, with some general results for $p$-local geometries, can be found in Benson and Smith \cite{bs}. In the present paper, it is understood that a $p$-local geometry for a group $G$ is a simplicial complex whose vertex stabilizers are suitably chosen maximal $p$-local subgroups of $G$.

\medskip
For many of the sporadic simple groups, especially when $p=2$ and $G$ has characteristic $2$-type, there is a $2$-local geometry $\Delta$, with vertex stabilizers maximal $2$-local subgroups of $G$, such that the Bouc complex $|\mathcal{B}_2(G)|$ is equal to the barycentric subdivision of the geometry. For the remaining sporadic groups the relationship is not clear anymore; in fact, $\Delta$ and $|\mathcal{B}_2(G)|$ are not always homotopy equivalent; see \cite{bs} for example.

\medskip
It is the purpose of this paper to give such a relationship for the Conway's third sporadic group $Co_3$. Section $2$ provides some notation and reviews standard results which will be used in proofs. In Section $3$ we describe in detail the $2$-local maximal ``parabolic" geometry $\Delta$ for $Co_3$. In Section $4$ we describe the relevant elements of the Bouc complex $|\mathcal{B}_2|$ of $Co_3$ in group theoretic and geometric terms. In Section $5$ we prove that the $2$-local geometry $\Delta$ is $G$-homotopy equivalent with a subcomplex $|\mathcal{\widehat B}_2|$ of the Bouc complex. In Section $6$ we prove that the fixed point set $\Delta ^z$ of a central involution $z$ is contractible. Then, using a result of Th\'{e}venaz \cite{th87}, we conclude that the reduced Lefschetz module $\tilde L ( \Delta)$ is projective relative to the collection $\mathfrak{X} = \mathcal{B}_2 \setminus \mathcal{\widehat B}_2$ of subgroups of $Co_3$.\\

\section{Notations, terminology and standard results}
\subsection{Subgroup complexes}Let $G$ be a finite group and $p$ a fixed prime which divides the order of $G$. A $p$-subgroup $R$ of $G$ is called $p$-radical if $R = O_p(N_G(R))$. The
$p$-subgroup $R$ is $p$-centric if $Z(R)$ is a Sylow $p$-subgroup of $C_G(R)$.

\medskip
A collection is a family of subgroups of $G$ which is closed under $G$-conjugation and it is ordered by inclusion; hence it is a $G$-poset. The following two collections are standard in the literature: \begin{align*}
\mathcal{A}_p (G) &= \lbrace E \mid E \; \text{nontrivial elementary abelian p-subgroup of} \; G  \rbrace,\\
\mathcal{B}_p (G) &= \lbrace R \mid R \; \text{nontrivial p-radical subgroup of} \; G \rbrace.
\end{align*}

In what follows, if $\mathcal{C}$ is a collection of subgroups of $G$ then $|\mathcal{C}|$ will denote the associated simplicial complex, with simplices equal to the flags or chains of subgroups in the poset $\mathcal{C}$. The complex $|\mathcal{A}_p(G)|$ is known as the Quillen complex and $|\mathcal{B}_p(G)|$ is known as the Bouc complex. The two complexes are $G$-homotopy
equivalent.

\medskip
Define the following subcollection of the Bouc collection:
$$\mathcal{\widehat B}_p(G) = \lbrace U \in \mathcal{B}_p(G) \mid Z(U) \cap Z(S) \not= 1, \; \text{for some} \; S \in \text{Syl}_p(G) \rbrace .$$
We shall call $|\mathcal{\widehat B}_p(G)|$ the {\it distinguished Bouc complex}.

\medskip
Let $\mathcal{E}_p(G)$ denote the collection of nontrivial elementary abelian $p$-subgroups of $G$ whose elements lie in the smallest set which contains elements of order $p$ in the center of a Sylow $p$-subgroup, is closed under conjugation in $G$ and is closed under taking products of commuting elements. This collection was introduced by Benson \cite{ben94} in order to study the mod-$2$ cohomology of $Co_3$.

\medskip
\subsection{Homotopy techniques}

Let us assume that $\Delta$ is a simplicial complex. We say that $G$ acts admissibly on $\Delta$ if whenever $g \in G$ fixes a simplex $\sigma \in \Delta$ then $g$ fixes every face of $\sigma$. If $v$ is a vertex of $\Delta$, then
$\text{Star}_{\Delta}(v)$ is the collection of simplices containing it,
and the residue of $v$, denoted $\text{Res}_{\Delta}(v)$ is the
subcomplex of next lower dimension, obtained by deleting $v$ from each
simplex of $\text{Star}_{\Delta}(v)$. Then $\text{Star}_{\Delta}(v)= v *
\text{Res}_{\Delta}(v)$, the simplicial join.

\medskip
For $\sigma \in \Delta$, let $G_{\sigma}$ be the stabilizer of $\sigma$
under this action and let $U_{\sigma}$ denote the kernel of the action of
$G_{\sigma}$ on the residue of $\sigma$. The following is inspired from a standard property of the finite groups of Lie type acting on their buildings:

\medskip
\textsc{The Borel-Tits property (BT)}: {\it Let $G$ be a group which acts admissibly on the simplicial complex $\Delta$. For each nontrivial
$p$-subgroup $U$ of $G$, there exists a simplex $\sigma$ of $\Delta$ such
that $N_G (U) \leq G_{\sigma}$}.

\begin{ntn}
For $P \leq G$, denote by $\Delta ^P$ the subcomplex of $\Delta$ fixed by
$P$.
\end{ntn}

The following two lemmas are standard results which will be used in the proofs. They are given here in the form used by Ryba, Smith and Yoshiara in \cite{rsy}.

\begin{lem}\cite[Lemma 2.1]{rsy} For a vertex $v$ of a simplicial
complex $\Delta$, if the residue $Res_{\Delta}(v)$ is contractible, then
$\Delta$ is homotopy equivalent to $\Delta \setminus Star(v)$.
\end{lem}

\begin{lem}\cite[Section 2]{rsy} Let $\Sigma \in \Delta$ be a simplex
of maximal dimension with $\sigma$ as a face. Assume that $\Sigma$ is the
only simplex of maximal dimension with $\sigma$ as a face. Then the process
of removing $\Sigma$ from $\Delta$, by collapsing $\Sigma$ down onto its
faces other than $\sigma$, is a homotopy equivalence.\\
\end{lem}

\section{The ${\mathbf 2}$-local geometry $\mathbf{\Delta}$ of $\mathbf{Co_3}$}

The sporadic simple group $G=Co_3$ has two classes of involutions \cite{fin73}, which we shall denote $2A$ (the central involutions) whose centralizer is $2^.S_6(2)$, and $2B$ (the non-central involutions) whose centralizer is $2 \times M_{12}$; using the Atlas \cite{atlas} notation. Note that central means that the elements are central in some Sylow $2$-subgroup of $G$. Given any two commuting central involutions, their product is a central involution. The product of a central involution and a non-central involution (if an involution) is a non-central involution. The product of two commuting non-central involutions can be either central or non-central.

\medskip
Let $\Delta$ denote the $2$-local maximal ``parabolic"
geometry of $G$. This geometry was first mentioned in \cite{rst}; further details can be found in \cite[Section 8.13]{bs}. This is a
rank $3$ geometry, whose objects will be denoted points $\mathcal{P}$,  lines $\mathcal{L}$, and $\mathcal{M}$-spaces. The objects of this geometry correspond to pure central elementary abelian subgroups of $G$; this means that these subgroups contain central involutions only. The points correspond to rank one subgroups $2$, lines to $2^2$ and $\mathcal{M}$-spaces to $2^4$. Incidence is given by
containment. $\Delta$ can also be regarded as a simplicial complex
of dimension two with three types of vertices.

\begin{rem}There is also a class of pure central elementary abelian subgroups $2^3$ in $G$, which do not correspond to objects in $\Delta$; each subgroup $2^3$ is contained in a
unique subgroup $2^4$; see \cite[Section 8.13]{bs}. These $2^3$ subgroups will be referred to as planes.
\end{rem}

The group $G$
acts (faithfully) flag-transitively on the geometry. The stabilizers of
the three types of objects are:

\medskip
\hspace*{1cm}\begin{tabular}{lll}
$G_p \simeq 2^.S_6(2)$ && for a point $p \in \mathcal{P}$;\\
$G_L \simeq 2^{2+6}3(S_3 \times S_3)$ && for a line $L \in \mathcal{L}$;\\
$G_M \simeq 2^{4.}L_4(2)$ && for an $\mathcal{M}$-space $M$.
\end{tabular}

\medskip
The flag stabilizers can be easily determined and they are:

\medskip
\hspace*{1cm}\begin{tabular}{llll}
$G_{pL} \simeq 2^{2+6}(S_3 \times S_3)$, &&\phantom{aaaaaaaa}&$G_{pM} \simeq 2^{1+6}_+L_3(2)$,\\
$G_{LM} \simeq  2^{2+6}(S_3 \times S_3)$, &&&$G_{pLM} \simeq \left [2^9\right ].S_3=\left [2^{10}3\right]$.
\end{tabular}

\medskip
The geometry $\Delta$ is pure $2$-local, in the sense of \cite{sy},
that is, all the simplex stabilizers have nontrivial normal $2$-subgroups. All vertex stabilizers are maximal $2$-local
subgroups. $\Delta$ has the following diagram:

\begin{center}
\begin{picture}(1000,50)(0,0)
\put(110,26){$\circ$}
\put(112,15){\scriptsize{$\mathcal{P}$}}
\put(115,29){\line(1,0){30}}
\put(145,26){$\circ$}
\put(152,15){\scriptsize{$\mathcal{L}$}}
\put(150,28){\line(1,0){30}}
\put(150,30){\line(1,0){30}}
\put(148,-3){\line(0,1){30}}
\put(132,-11){$\qed$}
\put(180,25){$\circ$}
\put(181,15){\scriptsize{$\mathcal{M}$}}
\end{picture}

\vspace*{1cm}
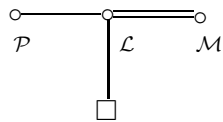
{\bf Figure 1:} {\it Diagram for the $2$-local geometry of $Co_3$}
\end{center}

\begin{ntn} For $X \in \lbrace \mathcal{P}, \mathcal{L}, \mathcal{M}
\rbrace$ and $F$ a flag of $\Delta$, we will denote by $X_F$ the
collection of all objects in $X$ incident with $F$.
\end{ntn}

\begin{ntn} For $p \in \mathcal{P}$ let
$$p^{\perp} = \lbrace q \in
\mathcal{P} \mid q \; \text{and} \; p \; \text{are incident with some
common line} \rbrace.$$
\end{ntn}

In what follows it will be useful to regard $\Delta$ as a point-line
geometry, that is, the lines and the $\mathcal{M}$-spaces are
identified with the subsets of points they are incident with. We give below some of the properties of $\Delta$:
\begin{itemize}
\item[$(\Delta 1).$] For $p \in \mathcal{P}$ the residue $\text{Res}_{\Delta}p = (
\mathcal{L}
_p, \mathcal{M}_p)$ is the collection of all lines and $\mathcal{M}$-spaces
containing $p$. Then $\text{Res}_{\Delta}p$ is
isomorphic with the geometry of isotropic lines and isotropic planes of a polar
space of type $C_3$ over the field $\mathbb{F}_2$. Note there are
$315$ lines and $135$ $\mathcal{M}$-spaces on a given point $p \in
\mathcal{P}$.
\item[$(\Delta 2).$]  For $M \in \mathcal{M}$, $\text{Res}_{\Delta} M = (\mathcal{P}_M, \mathcal{L}_M)$, the collection of points and
lines lying in $M$ is the truncation to points and lines of a projective
space $PG(3,2)$. It is immediate that there are $15$ points and $35$ lines
incident with $M$.
\item[$(\Delta 3).$]  For $L \in \mathcal{L}$, $\text{Res}_{\Delta}L = (\mathcal{P}_L, \mathcal{M}_L)$ is a digon, a complete bipartite graph on the three points and the three $\mathcal{M}$-spaces incident with $L$.
\item[$(\Delta 4).$]  Given a point $p$ and a line $L$ with $p^{\perp} \cap L \not = \emptyset$, either $p^{\perp} \cap L$ is a single point or all of $L$.
\item[$(\Delta 5).$]  If $(p, M) \in \mathcal{P} \times \mathcal{M}$, either $p \in M$, or $p^{\perp}
\cap M$ is at most a line.
\item[$(\Delta 6).$]  Two distinct $\mathcal{M}$-spaces intersect in at most one line.
\end{itemize}

The properties $(\Delta1)-(\Delta3)$ can be read from the diagram. For completeness we will provide arguments for $(\Delta4)-(\Delta6)$.\\
In order to prove $(\Delta4)$, consider $(p,L) \in \mathcal{P} \times \mathcal{L}$ with $p^{\perp} \cap L \not = \emptyset$. But $p^{\perp} \cap L$ cannot consist of two points only, because if the involution corresponding to the point $p$ commutes with two of the involutions in the $2^2$ subgroup corresponding to $L$, then it will commute with their product also. Thus $p^{\perp}\cap L$ can be a point or the entire line. The latter case occurs if $p$ and $L$ lie in the same $\mathcal{M}$-space.\\
Next we prove $(\Delta5)$. Let us assume that $p \not \in M$. It is easy to see that $p^{\perp}\cap M$ cannot be the entire $M$, since there are no projective subspaces of rank $4$ or greater in $\Delta$; also, $p^{\perp} \cap M$ cannot be a plane, by Remark $3.1$.\\
Finally, $(\Delta6)$ is a direct consequence of Remark $3.1$.

\begin{rem}
The Benson collection $\mathcal{E}_2=\mathcal{E}_2(G)$ contains
four conjugacy classes of subgroups, of orders $2, 2^2, 2^3, 2^4$; see \cite{ben94}. Let
$\mathcal{E}_2^-$ be the subcollection of $\mathcal{E}_2$ with the subgroups
of the form $2^3$ removed.  It is a direct consequence of Lemma $2.2$ and Remark $3.1$ that $|\mathcal{E}_2|$ and $|\mathcal{E}_2^-|$ are $G$-homotopy equivalent. Note that $|\mathcal{E}_2^-| = \Delta$. Simplices in $\Delta$ correspond to flags, or chains of subgroups in $\mathcal{E}_2^-$.
\end{rem}

\begin{rem}
Let $\tilde \Delta$ denote the involution geometry of $G$, the
point-line geometry whose points are the involutions of $G$. Two
involutions are collinear if they commute. Then $\Delta$ can be
identified with a (geometric) subspace of $\tilde \Delta$ whose
points
correspond to the central involutions, that is every exterior line of $\tilde
\Delta$ intersects $\Delta$ at a point or the empty set. Therefore,
we can regard $\Delta$ as a ``central involution geometry".\\
\end{rem}

\section{A geometric description of the distinguished Bouc complex of $\mathbf{Co_3}$}

In this section $G=Co_3$ and $\Delta$ represents its $2$-local geometry; further let $\mathcal{B} _2 = \mathcal{B} _2(Co_3)$ and $\mathcal{\widehat B} _2 = \mathcal {\widehat B} _2(Co_3)$. The Bouc complex of $Co_3$ was determined in \cite{an99}, from which we
reproduce the following table:\\

\begin{center}
\begin{tabular}{l l l l l}
{\bf Type} \phantom{aaaa}& {\bf Name} \phantom{aaaa}& $\mathbf{R}$ \phantom{aaaaaa}& $\mathbf{Z(R)}$ \phantom{aaaa}&
$\mathbf{N_G(R)}$\\
$p$&$R_p$ & $2$ & $A^1$ & $ 2^.S_6(2)$\\
$-$&$R_2$ & $2$ & $B^1$ & $2 \times M_{12}$\\
$-$&$R_3$ & $2^2$ & $B^3$ & $A_4 \times S_5$\\
$-$&$R_4$ & $2^3$ & $B^7$ & $(2^3 \times S_3)F^3_7$\\
$M$&$R_M$ & $2^4$ & $A^{15}$ & $2^4 A_8$\\
$p\square$&$R_{p \square}$ & $2^{1+5}$ & $A^1$ & $2^{1+5}S_6$\\
$pM$&$R_{pM}$ & $2^{1+6}_+$ & $A^1$ & $2^{1+6}_+ L_3(2)$\\
$M\square$&$R_{M \square}$ & $2^{3+4}$ & $A^7$ & $2^{3+4} L_3(2)$\\
$\mathbb{L}$&$R_L$ & $2^{2+6}$ & $A^3$ & $2^{2+6}3.(S_3 \times S_3)$\\
$pML$&$R_{pML}$ & $2^4.2^{2+3}$ & $A^1$ & $2^4.2^{2+3}.S_3$\\
$ML \square$&$R_{ML\square}$ & $2^4.2^{2+3}$ & $A^1$ &
$2^4.2^{2+3}.S_3$\\
$pM \square$&$R_{pM \square}$ & $2^4.2^{1+4}_+$ & $A^1$ & $2^4.2^{1+4}_+.S_3$\\
$pL \square$&$R_{pL\square}$ & $\left [2^9\right ]$ & $A^1$ & $\left [2^9\right ].S_3$\\
$pML \square$&$R_{pML\square}$ & $\left [2^{10}\right ]$ & $A^1$ & $\left [2^{10}\right ]$
\end{tabular}

\vspace*{.5cm}
{\bf Table 1:} {\it Representatives for conjugacy classes of radical $2$-subgroups of $Co_3$}
\end{center}

\bigskip
For $v \in \lbrace \mathcal{P}, \mathcal{L}, \mathcal{M} \rbrace$, the notation $R_v$ is suggested by the complex $\Delta$, specifically $R_v = O_2(G_v)$. The subgroup $R_L$ is not equal to a line $L$ in $\Delta$ (the pure central subgroups $2^2$ are not $2$-radical), but contains what we call a ``line structure" $\mathbb{L}$; see below for details. Recall that $\text{Res}_{\Delta}v$, for $v \in \lbrace \mathcal{P}, \mathcal{M} \rbrace$ is a ``truncation", obtained from a geometry of type $A_3$ or $C_3$ by ignoring one of the usual vertex types. We let $\square$ denote a vertex of that type {\it in the residue}; but note that $\square$ does not correspond to a vertex in the {\it full} geometry $\Delta$. In $\text{Res}_{\Delta}p, \; \square$ stands for a structure of $15$ lines and $15$ $\mathcal{M}$-spaces which form a generalized quadrangle $GQ(2,2)$. In $\text{Res}_{\Delta}M, \; \square$ can be identified with a plane of the projective space $PG(3,2)$. Groups involving $\square$, such as $R_{p \square}$ or $R_{M \square}$ can be constructed as the inverse images under the quotient maps $2^. S_6(2) \twoheadrightarrow S_6(2)$ or $2^{4.}L_4(2) \twoheadrightarrow L_4(2)$ of unipotent radicals of corresponding parabolics. The groups $R_{pM}$ and $R_{pML}$ can be constructed in the same manner, or as $O_2(G_{pM})$ and $O_2(G_{pML})$.

\medskip
Among the $2$-radical subgroups of $G$ all but the first four conjugacy classes of groups are $2$-centric. Also it is easy to see that $\mathcal{\widehat B}_2$ equals the collection obtained by removing from $\mathcal{B}_2$ the subgroups in the conjugacy classes of $\lbrace R_2, R_3, R_4 \rbrace$.

\medskip
As mentioned in the previous section, the standard $2$-local geometry
$\Delta$ can be regarded as a pure central involution geometry. In this section we will describe the subcomplex $|\mathcal{\widehat B}_2|$ in more geometric terms, using the properties of the central involutions of $G$.

\medskip
We start with a weaker version of the Borel-Tits property mentioned in Section $2$:

\medskip
(\textsc{BT})$^c$: {\it For each non-trivial $2$-subgroup $U$ of $G$ such that $Z(U)$
contains a central involution, there exists a simplex $\sigma$ of $\Delta$ such
that $N_G(U) \leq G_{\sigma}$}.

\medskip
\begin{prop} The pair $(G, \Delta)$ satisfies the property
\rm{(}\textsc{BT}\rm{)}$^c$.
\end{prop}

\begin{proof} Let $U \leq G $ be a non-trivial $2$-subgroup such that
$Z(U)$ contains a central involution. Let $H$ denote the set of central involutions in $Z(U)$. Recall that $H$ is closed under taking
products. Thus $H$ is an elementary abelian subgroup
of $U$ of rank at most $4$. Then $N_G(U) \leq N_G(H)$. If $H \simeq 2^3$ then $N_G(H)
\leq N_G(2^4)$. Otherwise $N_G(H) \in \lbrace G_p, G_L, G_M  \rbrace$ and the
conclusion follows.
\end{proof}

\begin{rem} Note that the pair $(G, \mathcal{\widehat B}_2)$ satisfies the property (\textsc{BT})$^c$. To see this let $U$ be a $2$-subgroup of $G$ with the property that $Z(U)$
contains a central involution. It is a direct consequence (of Proposition $4.1$ and the fact that the stabilizers of vertices in $\Delta$ are normalizers of groups $R$ which lie in $\widehat{\mathcal B}_2$) that there exists $R \in \widehat{\mathcal B}_2 (G)$ with $N_G(U) \leq N_G(R)$.
\end{rem}

In what follows we shall describe each of the subgroups of $\mathcal{\widehat B}_2$ in geometric terms; this means that we describe the set of central involutions contained in each distinguished radical $2$-subgroup. Recall a central involution is one lying in the center of a Sylow $2$-subgroup of $G$.

\begin{prop}The Sylow $2$-subgroup $R_{pML\square}$ of $G$ contains $55$ central involutions. These points of the geometry $\Delta$ which are in $R_{pML\square}$ lie either in a collection of three $\mathcal{M}$-spaces on a common line (which we refer to as a ``line structure" $\mathbb{L}$, consisting of $39$ points), or in a collection of $31$ points on $15$ lines, all containing a common cone point. These $31$ points also form $15$ planes, all containing the cone point, and these $15$ lines and $15$ planes can be thought of as forming a generalized quadrangle $GQ(2,2)$. The overlap $\mathbb{L} \cap GQ(2,2)$ consists of $15$ points, forming $3$ planes (one from each $\mathcal{M}$-space, and including their common line) and containing $7$ of the $15$ lines of $GQ(2,2)$.
\end{prop}

\begin{proof} A computation using GAP \cite{gap} verified the existence of the $55$ central involutions. The Sylow $2$-subgroup $R_{pML\square}$ of $G$ can be described as an extension $2^4.U_4$, where $U_4$ is the group of upper triangular matrices in $L_4(2)$. Clearly $M=2^4$ contains $15$ central involutions; we are looking for $40$ others. Let $z \in R_{pML\square}$ be a central involution such that $z \not \in M$, and denote by $\overline{z}$ the image of $z$ in the quotient group $U_4$. It can be shown that all of the involutions in $U_4$ lie either in the elementary abelian subgroup

{\small $$\left\{\left(
  \begin{array}{cccc}
    1 & 0 & * & * \\
    0 & 1 & * & * \\
    0 & 0 & 1 & 0 \\
    0 & 0 & 0 & 1 \\
  \end{array}
\right)\right\}$$}

of rank four, or in the extraspecial group

{\small $$ 2_+^{1+4} \; = \;\left\{\left(
  \begin{array}{cccc}
    1 & * & * & * \\
    0 & 1 & 0 & * \\
    0 & 0 & 1 & * \\
    0 & 0 & 0 & 1 \\
  \end{array}
\right)\right\}.$$}

Write $\overline{z}=I+N$ where $I$ is the identity matrix and $N$ is nilpotent. In fact, $N^2=0$ since $\overline{z}$ has order two.

\medskip
The geometric condition that $p^{\perp} \cap M$ is at most a line, see Property $(\Delta5)$, Section $3$, implies that the set of elements of $M$ fixed by the action of the matrix $\overline{z}$ is at most $2^2$. Thus $\overline{z}$ cannot be a transvection, which would fix a plane $2^3$. Since $\overline{z}v=v$ if and only if $Nv=0$, we must have rank$(N)\not= 1$. There are precisely ten upper triangular $4 \times 4$ matrices satisfying rank$(N)=2$ and $N^2=0$, namely:

{\small $$N_1 = \left(
  \begin{array}{cccc}
    0 & 0 & 1 & 0 \\
    0 & 0 & 0 & 1 \\
    0 & 0 & 0 & 0 \\
    0 & 0 & 0 & 0 \\
  \end{array}
\right)\;
N_2 = \left(
  \begin{array}{cccc}
    0 & 0 & 1 & 1 \\
    0 & 0 & 1 & 0 \\
    0 & 0 & 0 & 0 \\
    0 & 0 & 0 & 0 \\
  \end{array}
\right)\;
N_3 = \left(
  \begin{array}{cccc}
    0 & 0 & 0 & 1 \\
    0 & 0 & 1 & 1 \\
    0 & 0 & 0 & 0 \\
    0 & 0 & 0 & 0 \\
  \end{array}
\right)$$}

{\small $$N_4 = \left(
  \begin{array}{cccc}
    0 & 0 & 0 & 1 \\
    0 & 0 & 1 & 0 \\
    0 & 0 & 0 & 0 \\
    0 & 0 & 0 & 0 \\
  \end{array}
\right)\quad
N_5 = \left(
  \begin{array}{cccc}
    0 & 0 & 1 & 1 \\
    0 & 0 & 0 & 1 \\
    0 & 0 & 0 & 0 \\
    0 & 0 & 0 & 0 \\
  \end{array}
\right)\quad
N_6 = \left(
  \begin{array}{cccc}
    0 & 0 & 1 & 0 \\
    0 & 0 & 1 & 1 \\
    0 & 0 & 0 & 0 \\
    0 & 0 & 0 & 0 \\
  \end{array}
\right)$$}

{\small $$N_7 = \left(
  \begin{array}{cccc}
    0 & 1 & 0 & 0 \\
    0 & 0 & 0 & 0 \\
    0 & 0 & 0 & 1 \\
    0 & 0 & 0 & 0 \\
  \end{array}
\right)\quad
N_8 = \left(
  \begin{array}{cccc}
    0 & 1 & 0 & 1 \\
    0 & 0 & 0 & 0 \\
    0 & 0 & 0 & 1 \\
    0 & 0 & 0 & 0 \\
  \end{array}
\right)\quad
N_9 = \left(
  \begin{array}{cccc}
    0 & 1 & 1 & 0 \\
    0 & 0 & 0 & 1 \\
    0 & 0 & 0 & 1 \\
    0 & 0 & 0 & 0 \\
  \end{array}
\right)$$}

{\small $$N_{10} = \left(
  \begin{array}{cccc}
    0 & 1 & 1 & 1 \\
    0 & 0 & 0 & 1 \\
    0 & 0 & 0 & 1 \\
    0 & 0 & 0 & 0 \\
  \end{array}
\right)$$}

Denote by $\overline{z}_i, 1 \leq i \leq 10$, the corresponding unipotent matrix $\overline{z}_i=I + N_i$ (this is an abuse of notation, since we have not yet defined the element $z_i$ in $R_{pML\square}$). Note that $\overline{z}_3 = \overline{z}_1 \cdot \overline{z}_2 = \overline{z}_2 \cdot \overline{z}_1$, and we also have the commuting products $\overline{z}_6 = \overline{z}_4 \cdot \overline{z}_5$, $\overline{z}_9  = \overline{z}_1 \cdot \overline{z}_8 = \overline{z}_5 \cdot \overline{z}_7$ and $\overline{z}_{10} = \overline{z}_1 \cdot \overline{z}_7 = \overline{z}_5 \cdot \overline{z}_8$. Denote $M=2^4= \langle a_1, a_2, a_3, a_4 \rangle$ as well as $p=2=\langle a_1 \rangle$, $L = 2^2 = \langle a_1, a_2 \rangle$, and $\square = 2^3 = \langle a_1, a_2, a_3 \rangle$. We see that $\overline{z}_1, \overline{z}_2, \ldots \overline{z}_6$ each fix pointwise the line $\langle a_1, a_2 \rangle$, $\overline{z}_7$ and $\overline{z}_8$ fix the line $\langle a_1, a_3 \rangle$, and $\overline{z}_9$ and $\overline{z}_{10}$ fix the line $\langle a_1, a_2 a_3 \rangle$.

\medskip
If $z$ and $z'$, not in $M$, are two distinct central involutions in $R_{pML\square}$ with the same image $\overline{z}=\overline{z'}$ in $U_4$, then $z \cdot z' \in M$, implying $z \cdot z' \in z^{\perp} \cap M$. Since $z^{\perp} \cap M$ is at most one line containing $3$ points, for each point $z$ there are at most three other points $z'$ with $\overline{z} = \overline{z'}$. We are counting the four elements in a coset $z \cdot 2^2$. Thus we have found at most $40= 4 \times 10$ central involutions in $z \in R_{pML\square}$ with $z \not \in M$. Combining this with the fact that $R_{pML\square}$ contains $55$ central involutions (the aforementioned GAP computation) implies that we have described exactly $40$ such central involutions not in $M$.

\medskip
Let $z_i, 1 \leq i \leq 10$, denote a representative central involution in $R_{pML\square}$ with image the matrix $\overline{z}_i \in U_4$.

\medskip
Note that $M_1 = M = \langle a_1, a_2, a_3, a_4 \rangle$, $M_2 = \langle a_1, a_2, z_1, z_2 \rangle$, and $M_3 = \langle a_1, a_2, z_4, z_5 \rangle$ are three $\mathcal{M}$-spaces on the common line $L= \langle a_1, a_2 \rangle$. These form the line structure $\mathbb{L}$, containing $39$ points.

\medskip
Let us now describe the generalized quadrangle $GQ(2,2)$. The fifteen lines, each spanned by $p=\langle a_1 \rangle$ and one element of the following set $\lbrace a_2$, $a_3$, $a_2a_3$, $z_1$, $a_2z_1$, $z_5$, $a_2z_5$, $z_7$, $a_3z_7$, $z_8$, $a_3z_8$, $z_9$, $a_2a_3z_9$, $z_{10}$, $a_2a_3z_{10} \rbrace$, are the $15$ lines on a common cone point $p$. The fifteen planes $\square _i, 1 \leq i \leq 15$ of the generalized quadrangle $GQ(2,2)$ are spanned by $p= \langle a_1 \rangle$ and one of $\lbrace \langle a_2,a_3 \rangle$, $\langle a_2,z_1 \rangle$, $\langle a_2,z_5 \rangle$, $\langle a_3,z_7 \rangle$, $\langle a_3,z_8 \rangle$, $\langle a_2a_3,z_9 \rangle$, $\langle a_2a_3,z_{10} \rangle$, $\langle z_1,z_7 \rangle$, $\langle z_1,z_8 \rangle$, $\langle z_5,z_7 \rangle$, $\langle z_5,z_8 \rangle$, $\langle a_2z_1,a_3z_7 \rangle$, $\langle a_2z_1,a_3z_8 \rangle$, $\langle a_2z_5, a_3 z_7 \rangle$, $\langle a_2z_5,a_3z_8 \rangle\rbrace$.

\medskip
The overlap $\mathbb{L} \cap GQ(2,2)$ equals $\square _1 \cup \square _2 \cup \square _3$ where $\square _1 = \square = \langle a_1, a_2, a_3 \rangle$, $\square _2 = \langle a_1, a_2, z_1 \rangle$, and $\square _3 = \langle a_1, a_2, z_5 \rangle$. Note that $\square _i \subseteq M_i$ for $1 \leq i \leq 3$. Also $M_1, M_2, M_3, \square _4, \square _5, \ldots \square _{15}$ are the maximal pure central elementary abelian subgroups of $R_{pML\square}$.
\end{proof}

\bigskip
We now describe the other $2$-radical subgroups of $G$, from both geometric and group theoretic points of view.

\medskip
${\mathbf{1)}}\;R_p \simeq 2$, generated by a central involution (which the Atlas \cite{atlas} lists as the conjugacy class $2A$), is a point of the geometry $\Delta$.

\medskip
${\mathbf{2)}}\;R_M \simeq 2^4$, a pure central elementary abelian $2$-group of rank four, is an $\mathcal{M}$-space of the geometry $\Delta$. It contains $15$ points, or central involutions, as well as $35$ lines and $15$ planes.

\medskip
${\mathbf{3)}}\;R_ {p \square} \simeq 2^{1+5}$ is defined as a subgroup of $G_p = 2^.Sp_6(2)$ as the inverse image of the unipotent radical of the parabolic subgroup $2^5.S_6$ in the quotient group $Sp_6(2)$. Under the quotient map, a central involution of type $2A$ lying in $2^.Sp_6(2)$ maps either to the identity or to an involution in the conjugacy class $2B$ in $Sp_6(2)$, using the Atlas \cite{atlas} notation. A computation using GAP \cite{gap} shows that $R_{p\square}$ contains $31$ central involutions (there are $15$ involutions of type $2B$ in $2^5 \subseteq Sp_6(2)$), and that these $31$ points lie on the $15$ lines through a cone point, the involution in the center $Z(R_{p\square})$, and form the $15$ planes of the generalized quadrangle $GQ(2,2)$. As a subgroup of a Sylow $2$-subgroup $R_{pML\square}=2^4.U_4$, $R_{p\square}$ is generated by $\langle a_1, a_2, a_3, z_1, z_5, z_7 \rangle$, with image $\langle \overline{z}_1, \overline{z}_5, \overline{z}_7 \rangle$ an elementary abelian subgroup of $U_4$ of rank three. Thus $R_{p \square} \simeq 2^3.2^3 \subseteq 2^4.U_4$. $R_{p \square}$ does not contain any $\mathcal{M}$-spaces. There are a total of $75$ lines in $R_{p \square}$, the $15$ lines on the cone point and $60$ others, contained in the $15$ planes but not through the cone point.

\medskip
${\mathbf{4)}}\;R_{pM} \simeq 2^{1+6}_+$ can be defined as a subgroup of $G_p=2^.Sp_6(2)$ as the inverse image of the unipotent radical of the parabolic subgroup $2^6.L_3(2)$ in the quotient group $Sp_6(2)$. However, $R_{pM}$ can also be defined as a subgroup of $G_M=2^4.L_4(2)$ as the inverse image of the unipotent radical of a parabolic subgroup $2^3.L_3(2)$ in the quotient group $L_4(2)$. As a subgroup of a Sylow $2$-subgroup $R_{pML\square}$, $R_{pM} \simeq 2^4.2^3 \subseteq 2^4.U_4$, with quotient group $2^3 \subseteq U_4$ the group of matrices

{\small $$\left\{ \left(
  \begin{array}{cccc}
    1 & * & * & * \\
    0 & 1 & 0 & 0 \\
    0 & 0 & 1 & 0 \\
    0 & 0 & 0 & 1 \\
  \end{array}
\right)\right\}.$$}

Thus $R_{pM}$ contains only the $15$ points of $M= \langle a_1, a_2, a_3, a_4 \rangle$. Of course, $R_M$ is a subgroup of $R_{pM}$. Clearly $R_{pM}$ does not contain any subgroup conjugate to $R_{p\square}$ since the latter contains $31$ points.

\medskip
${\mathbf{5)}}\;R_{M\square} \simeq 2^{3+4} \simeq 2^4.2^3 \subseteq 2^4.U_4$ can be defined as a subgroup of $G_M=2^4.L_4(2)$ as the inverse image of the unipotent radical of a parabolic subgroup $2^3.L_3(2)$, stabilizing a plane, in the quotient group $L_4(2)$. The group $2^3 \subseteq U_4$ is the group of matrices

{\small $$\left\{ \left(
  \begin{array}{cccc}
    1 & 0 & 0 & * \\
    0 & 1 & 0 & * \\
    0 & 0 & 1 & * \\
    0 & 0 & 0 & 1 \\
  \end{array}
\right)\right\},$$}

and thus $R_{M\square}$ contains only the $15$ points (and the $\mathcal{M}$-space) of $M=\langle a_1, a_2, a_3, a_4 \rangle$. Note $R_{M\square}$ is not isomorphic to $R_{pM}$ since their centers have different orders; $Z(R_{M\square})=2^3=\langle a_1,a_2,a_3 \rangle$ and $Z(R_{pM})=2=\langle a_1 \rangle$. As above, $R_{M\square}$ does not contain any conjugate of $R_{p\square}$.

\medskip
${\mathbf{6)}}\;R_L \simeq 2^{2+6}$ equals $O_2(G_L)=O_2(N_G(2^2))$, and has center $Z(R_L)$ the line $L=2^2=\langle a_1, a_2 \rangle$. Note that $G_L=2^{2+6}3(S_3 \times S_3)$ with the first $S_3$ permuting the $3$ points of $L$ and the second $S_3$ permuting the three $\mathcal{M}$-spaces on $L$. Let $M$ be any of the three $\mathcal{M}$-spaces on $L$; then $G_{LM}=G_L \cap G_M= 2^{2+6}3(S_3 \times 2)$, which is isomorphic to $2^{2+6}(S_3 \times S_3)$. Therefore $R_M \subseteq O_2(G_{LM}) = O_2(G_L) = R_L$. Thus $R_L$ contains all three of the $\mathcal{M}$-spaces on $L$, which is the line structure $\mathbb{L}$ consisting of $39$ points. Also, $R_L$ is the subgroup of $G_M=2^4.L_4(2)$ which equals the inverse image of the unipotent radical of the parabolic subgroup $2^4(S_3 \times S_3)$ in $L_4(2)$, stabilizing a line $L$. As a subgroup of a Sylow $2$-subgroup $R_{pML\square}$, $R_L \simeq 2^4.2^4 \subseteq 2^4.U_4$, with quotient group $2^4 \subseteq U_4$ the group of matrices

{\small $$\left\{ \left(
  \begin{array}{cccc}
    1 & 0 & * & * \\
    0 & 1 & * & * \\
    0 & 0 & 1 & 0 \\
    0 & 0 & 0 & 1 \\
  \end{array}
\right)\right\}.$$}

It follows that $R_L$ contains precisely $39$ central involutions. Also $L$ is generated by its central involutions, $L = \langle a_1, a_2, a_3, a_4, z_1, z_2, z_4, z_5 \rangle$. Further $R_L$ does not contain any conjugate of $R_{p\square}$, since if $R_{p\square} \subseteq R_L$ then at least $11$ of the $31$ points of $R_{p\square}$ would have to lie on one of the $\mathcal{M}$-spaces $M$ of $R_L$. Since $11 > 7$, these points would span $M$. But then $M \subseteq R_{p\square}$, a contradiction, since $R_{p\square}$ does not contain any $\mathcal{M}$-spaces.

\medskip
${\mathbf{7)}}\;R_{pML} \simeq 2^4.2^{2+3} \subseteq 2^4.U_4$, with quotient group $2^{2+3} \subseteq U_4$ the group of matrices

{\small $$\left\{ \left(
  \begin{array}{cccc}
    1 & * & * & * \\
    0 & 1 & * & * \\
    0 & 0 & 1 & 0 \\
    0 & 0 & 0 & 1 \\
  \end{array}
\right)\right\}.$$}

$R_{pML}$ has precisely $39$ points and three $\mathcal{M}$-spaces of the line structure $\mathbb{L}$ contained in $R_L \subseteq R_{pML}$.

\medskip
${\mathbf{8)}}\;R_{ML\square} \simeq 2^4.2^{2+3} \subseteq 2^4.U_4$, with quotient group $2^{2+3} \subseteq U_4$ the group of matrices

{\small $$\left\{ \left(
  \begin{array}{cccc}
    1 & 0 & * & * \\
    0 & 1 & * & * \\
    0 & 0 & 1 & * \\
    0 & 0 & 0 & 1 \\
  \end{array}
\right)\right\}.$$}

$R_{ML\square}$ contains precisely the $39$ points and the three $\mathcal{M}$-spaces of the line structure $\mathbb{L}$ in $R_L \subseteq R_{pML}$.

\medskip
${\mathbf{9)}}\;R_{pM\square} = 2^4.2^{1+4}_+\subseteq 2^4.U_4$ with quotient group $2^{1+4}_+ \subseteq U_4$ the group of matrices

{\small $$\left\{ \left(
  \begin{array}{cccc}
    1 & * & * & * \\
    0 & 1 & 0 & * \\
    0 & 0 & 1 & * \\
    0 & 0 & 0 & 1 \\
  \end{array}
\right)\right\}.$$}

This implies that $R_{pM\square}$ contains $39$ central involutions, with $31$ points in $R_{p\square} \subseteq R_{pM\square}$ and $15$ points in $R_M \subseteq R_{pM}$. These two sets intersect in a set of $7$ points, the plane $\square = \langle a_1, a_2, a_3 \rangle \subseteq M$.

\medskip
${\mathbf{10)}}\;R_{pL\square} = 2^4.2^{2+3}\subseteq 2^4.U_4$ contains all $55$ points of a Sylow $2$-subgroup $R_{pML\square}$ since $R_{p\square} \subseteq R_{pL\square}$ and $R_L \subseteq R_{pL\square}$. Using $z_7 \in R_{pL\square}$ and $R_L \subseteq R_{pL\square}$, it follows that the quotient group $2^{2+3} \subseteq U_4$ equals the group of matrices

{\small $$\left\{ \left(
  \begin{array}{cccc}
    1 & \epsilon & * & * \\
    0 & 1 & * & * \\
    0 & 0 & 1 & \epsilon \\
    0 & 0 & 0 & 1 \\
  \end{array}
\right)\right\},$$}

the group of upper triangular $4 \times 4$ unipotent matrices $(x_{ij})$ satisfying $x_{12}= x_{34}$.

\medskip
${\mathbf{11)}}\;R_{pML\square} = 2^4.U_4$ contains a unique conjugate of each of the groups $R_{pL\square}$, $R_{pM\square}$, $R_L$ and $R_{p\square}$. This can be seen for the groups $R_{pL\square}$, $R_L$, and $R_{p\square}$ by using their geometric structures, involving the line structure $\mathbb{L}$ and the generalized quadrangle $GQ(2,2)$, and the fact that these three groups are generated by their central involutions. For $R_{pM\square}$, a simple counting argument applies. If $R_{pM\square} \subseteq R_{pML\square}$ then $R_{pML\square} \subseteq N_G(R_{pM\square}) \simeq R_{pM\square} .S_3$, and the latter group contains three Sylow $2$-subgroups.\\

\section{The homotopy relation}

In this section we prove that the $2$-local geometry $\Delta$ is $G$-homotopy equivalent to the distinguished Bouc complex $|\widehat{\mathcal{B}}_2|$ of $G=Co_3$. The complex $|\widehat{\mathcal{B}}_2|$ can be equivariantly retracted to a subcomplex $\Delta_1$ formed from a subposet of $\widehat{\mathcal{B}}_2$. Also $\Delta_1$ can be retracted to a subcomplex $\Delta_0$ which is homeomorphic to $\Delta$.

\medskip
Let $\Delta_1$ denote the simplicial complex corresponding to the subposet of $\mathcal{B}_2$ with the three types of objects $\lbrace R_p, R_M, R_L \rbrace$. Let $\Delta _0$ denote the subcomplex of $\Delta_1$ obtained by removing those simplices corresponding to the chains $R_p \subseteq R_L$ and $R_p \subseteq R_M \subseteq R_L$ where the subgroup $R_p$ is not normal in the group $R_L$.

\begin{prop}
(a). $\Delta _0$ is homeomorphic to $\Delta$.\\
\hspace*{3.6cm}(b). $\Delta _0$ is $G$-homotopy equivalent to $\Delta _1$.
\end{prop}

\begin{proof} (a). The homeomorphism from $\Delta$ to $\Delta _0$ is induced by the correspondence $p \rightarrow R_p$ (recall $p=R_p=2$), $L \rightarrow R_L$ (recall $R_L=2^{2+6}=O_2(N_G(L))$ with $L=2^2=Z(R_L)$), and $M \rightarrow R_M$ (recall $M=R_M=2^4$). Next, we show that simplices in $\Delta$ correspond to simplices in $\Delta _0$. Clearly $p \subseteq M$ corresponds to $R_p \subseteq R_M$. Also, $R_p \trianglelefteq R_L$ if and only if $R_p \subseteq Z(R_L)$ if and only if $p \subseteq L$. Finally, $L \subseteq M$ if and only if $R_M \subseteq R_L$, and $p \subseteq L \subseteq M$ corresponds to $R_p \subseteq R_M \subseteq R_L$ with $R_p \trianglelefteq R_L$.\\
(b). We use Lemma $2.3$ to remove those simplices in $\Delta_1$, not lying in $\Delta _0$, corresponding to chains $R_p \subseteq R_L$ and $R_p \subseteq R_M \subseteq R_L$ where the subgroup $R_p$ is not normal in $R_L$, so that $R_p \not\subseteq Z(R_L) =L$. Of course $R_p$ does commute with elements of the center $Z(R_L)$, and they generate a group $2^3= \langle R_p, Z(R_L) \rangle$. This $2^3$ lies in a unique $2^4=M$, an $\mathcal{M}$-space satisfying $L \subseteq M$. Therefore the chain $R_p \subseteq R_L$ is a face of a unique chain $R_p \subseteq R_M \subseteq R_L$, and Lemma $2.3$ allows to cancel these two simplices.
\end{proof}

\bigskip
In this section we will prove the main theorem of the paper. The proof of this theorem will be the focus of all of this Section $5$.

\begin{ta}{\it The two complexes $\Delta$ and $| \widehat{\mathcal{B}}_2|$ are $G$-homotopy equivalent.}
\end{ta}

We will reduce $|\widehat{\mathcal{B}}_2|$ down to the subcomplex $\Delta_1$ in a number of steps; at each step we will perform a homotopy retraction, which will use one of the results given in Lemmas $2.2$ or $2.3$. We first apply Lemma $2.2$ to retract $|\widehat{\mathcal{B}}_2|$ to a subcomplex $|\widehat{\mathcal{B}}_2^I|$ by removing all simplices (flags, or chains of subgroups) which involve groups of type $R_{pM}$ or $R_{M\square}$. Then Lemma $2.2$ will allow us to retract $|\widehat{\mathcal{B}}_2^I|$ to a subcomplex $|\widehat{\mathcal{B}}_2^{II}|$ by removing all simplices which involve groups of type $R_{pML}$ or $R_{ML\square}$.

\medskip
The main part of the proof involves the use of Lemma $2.3$ to reduce $|\widehat{\mathcal{B}}_2^{II}|$ down to a subcomplex $\Delta_2$ (in a sequence of $35$ steps), which in turn can easily be retracted to $\Delta _1$. We will give in Table $2$ a list of all orbits of flags that lie in $|\widehat{\mathcal{B}}_2^{II}|$ but not in $\Delta _1$, and we proceed to cancel all but $44$ of them to arrive at $\Delta_2$. We then apply Lemma $2.3$ to remove these $44$ flags as well.

\begin{prop} $|\widehat{\mathcal{B}}_2|$ is $G$-homotopy equivalent to the subcomplex $|\widehat{\mathcal{B}}_2^I|$, where $\widehat{\mathcal{B}}_2^I$ is the subposet of $\widehat{\mathcal{B}}_2$ obtained by removing the groups of the form $R_{pM}$ and $R_{M\square}$.
\end{prop}

\begin{proof} There exists a unique $R_M$ in $R_{pM}$, and there are exactly $15$ central involutions in $R_{pM}$, all contained in $R_M$. It follows that any group $R$ in $\widehat{\mathcal{B}}_2$ which is incident with $R_{pM}$ (either $R \subseteq R_{pM}$ or $R_{pM} \subseteq R$) is also incident with $R_M$. In the simplicial complex $|\widehat{\mathcal{B}}_2|$, the residue of the vertex $R_{pM}$ is a cone on the vertex $R_M$, so is contractible. According to Lemma $2.2$, we can remove the vertices of the form $R_{pM}$ together with their stars.

\medskip
Similarly, there exist a unique $R_M$ in $R_{M\square}$, containing all $15$ central involutions in $R_{M\square}$. Any group $R$ in $\widehat{\mathcal{B}}_2$ which is incident with $R_{M\square}$ is also incident with $R_M$. In $|\widehat{\mathcal{B}}_2|$, the residue of the vertex $R_{M\square}$ is contractible, a cone on the vertex $R_M$. Another application of Lemma $2.2$ gives the result.
\end{proof}

\begin{prop} $|\widehat{\mathcal{B}}_2^I|$ is $G$-homotopy equivalent to a subcomplex $|\widehat{\mathcal{B}}_2^{II}|$, where $\widehat{\mathcal{B}}_2^{II}$ is the subposet of $\widehat{\mathcal{B}}_2^I$ obtained by removing the groups of the form $R_{pML}$ and $R_{ML\square}$.
\end{prop}

\begin{proof} There exists a unique $R_L$ in $R_{pML}$ which contains all $39$ central involutions in $R_{pML}$. In the poset $\widehat{\mathcal{B}}_2^I$, having already removed
the groups of type $R_{pM}$, any group $R$ in $\widehat{\mathcal{B}}_2^I$ which is incident with $R_{pML}$ is also incident with $R_L$. In the simplicial complex $|\widehat{\mathcal{B}}_2^I|$, the residue of the vertex $R_{pML}$ is contractible, a cone on $R_L$. Similarly, there is a unique $R_L$ in $R_{ML\square}$ which contains the $39$ central involutions of $R_{pML}$. Since we have removed the groups of type $R_{M\square}$, any group $R$ in $\widehat{\mathcal{B}}_2^I$ which is incident with $R_{ML\square}$ is also incident with $R_L$. In $|\widehat{\mathcal{B}}_2^I|$, the residue of $R_{ML\square}$ is a cone on $R_L$. Therefore two applications of Lemma $2.2$ prove the Proposition.
\end{proof}

\medskip
Let us remark on the importance of doing these retractions in a particular order. We must remove the groups of type $R_{pM}$ before removing the groups of type $R_{pML}$, since in the original poset $\widehat{\mathcal{B}}_2$, $R_{pM}$ is incident with $R_{pML}$ but not with $R_L$.

\begin{center}
\begin{small}
\begin{center}
\begin{tabular}{|l|l|l||l|l|l|}
\hline
&&&&&\\
{\bf Rank 1 flags}&&&{\bf Rank 3 flags}&& \\

$p\square$ & $2^{10} \cdot 3^2 \cdot 5$ & $1$ & $p \subseteq M \subseteq pM\square$ & $2^{10} \cdot 3$ & $1,6,8$\\

$pM\square$ & $2^{10} \cdot 3$ & $1$ & $p\subseteq M \subseteq pL\square$ & $2^{10} \cdot 3$ & $1,2,12,24$\\

$p L \square$ & $2^{10} \cdot 3$ & $1$ & $p \subseteq M \subseteq pML \square$ & $2^{10}$ & $1,2,4,8,8,16$ \\

$pML\square$ & $2^{10}$ & $1$ &  $p \subseteq {\mathbb L} \subseteq pL\square$ & $2^{10} \cdot 3$
& $1,2,12,24$ \\

&&& $p\subseteq {\mathbb L} \subseteq pML\square$ & $2^{10}$ &
$1,2,4,8,8,16$ \\

{\bf Rank 2 flags}&&& $p \subseteq p\square \subseteq pL\square$ & $2^{10} \cdot 3$
& $1,2,12,16$ \\

$p \subseteq p\square$ & $2^{10} \cdot 3^2 \cdot 5$ & $1,30$ &$p \subseteq p\square \subseteq pM \square$ & $2^{10} \cdot 3$ &
$1,6,24$ \\

$p \subseteq pM\square$ & $2^{10} \cdot 3$ & $1,6,8,24$ &$p \subseteq p\square \subseteq pML\square$ & $2^{10}$ & $1,2,4,8,16$\\

$p \subseteq pL\square$ & $2^{10} \cdot 3$ & $1,2,12,16,24$ &$p \subseteq pM\square \subseteq pML \square$ & $2^{10}$ & $1,2,4,8,8,16$ \\

$p \subseteq pML\square$ & $2^{10}$ & $1,2,4,8,8,16,16$ &$p \subseteq pL\square \subseteq pML \square$ & $2^{10}$ & $1,2,4,8,8,16,16$ \\

$M \subseteq pM\square$ & $2^{10} \cdot 3$ & $1$ &$M \subseteq {\mathbb L} \subseteq pL\square$ & $2^{10} \cdot 3$ & $3$ \\

$M \subseteq pL\square$ & $2^{10} \cdot 3$ & $3$ &$M \subseteq {\mathbb L} \subseteq pML \square$ & $2^{10}$ & $1,2$ \\

$M \subseteq pML\square$ & $2^{10}$ & $1,2$&$M \subseteq pM\square \subseteq pML\square$ & $2^{10}$ & $1$ \\

${\mathbb L} \subseteq p L \square$ & $2^{10} \cdot 3$ & $1$
&$M \subseteq pL\square \subseteq pML\square$ & $2^{10}$ & $1,2$ \\

${\mathbb L} \subseteq pML\square$ & $2^{10}$ & $1$ &${\mathbb L} \subseteq pL\square \subseteq pML
\square$ & $2^{10}$ & $1$ \\

$p\square \subseteq pM\square$ & $2^{10} \cdot 3$ & $1$  &$p\square \subseteq pL\square \subseteq pML \square$
& $2^{10}$ & $1$ \\

$p\square \subseteq pL\square$ & $2^{10} \cdot 3$ & $1$   &$p\square \subseteq pM\square \subseteq pML\square$ & $2^{10}$ & $1$ \\

$p\square \subseteq pML\square$ & $2^{10}$ & $1$ &&&\\

$pM\square \subseteq pML\square$ & $2^{10}$ & $1$ &{\bf Rank 4 flags}&&\\

$pL\square \subseteq pM L \square$ & $2^{10}$ & $1$  &$p \subseteq M \subseteq {\mathbb L} \subseteq p L \square$ & $2^{10} \cdot 3$ & $1,2,12,24$ \\

&&&$p \subseteq M \subseteq {\mathbb L} \subseteq pM L \square$   &$2^{10}$&$1,2,4,8,8,16$\\

&&&$p \subseteq M \subseteq pM \square \subseteq pM L \square$ &
$2^{10}$ & $1,2,4,8$\\

&&&$p \subseteq M \subseteq pL\square \subseteq pML\square$ &
$2^{10}$ & $1,2,4,8,8,16$\\

&&&$p \subseteq {\mathbb L} \subseteq pL \square \subseteq pM
L \square$ & $2^{10}$ & $1,2,4,8,8,16$\\

&&& $p \subseteq p \square \subseteq pM \square \subseteq pML
\square$ & $2^{10}$ & $1,2,4,8,16$\\

&&& $p \subseteq p \square \subseteq pL \square \subseteq pM
L \square$ & $2^{10}$ & $1,2,4,8,16$\\

&&&$M \subseteq {\mathbb L} \subseteq p L \square \subseteq pM
L \square$ & $2^{10}$ & $1,2$\\

&&&& &\\
&&&{\bf Rank 5 flags}& &\\

&&&$p \subseteq M \subseteq {\mathbb L} \subseteq p L \square
\subseteq pM L \square$ & $2^{10}$ & $1,2,4,8,8,16$\\
&&&&&\\
\hline
\end{tabular}
\end{center}
\end{small}

\vspace*{.3cm}
{\bf Table 2:} {\it Orbits of flags in $|\widehat{\mathcal{B}}_2^{II}|$ which do not lie in $\Delta_1$}\end{center}

\subsection{Flags in $|\mathcal{\widehat B}_2^{II}|$}Let $F: \; x_1 \subseteq x_2 \ldots \subseteq x_n$ be a flag in $|\mathcal{\widehat B}^{II}_2|$. Each group $x_i$ will have a type which is a subset of $\lbrace p, M, L, \square \rbrace$, but we sometimes abbreviate this by saying $x_i$ has ``type $i$". The normalizer of the group $x_n$ acts on the collection of groups of type $i$ contained in $x_n$; in many cases this action is not transitive. In the first column  of Table $2$ we give the type of a flag $F$, for those orbits of flags in $|\widehat{\mathcal{B}}_2^{II}|$ but not in $\Delta _1$ (so that each flag contains a group $x_i$ whose type involves $\square$). In the second column of Table $2$ we give the order of the normalizer of the group $x_n$ of $F$, and in the third column we give the number of elements in the orbits of the action of $N_G(x_n)$ on the collection of groups of type $1$ contained in $x_n$.

\medskip
To prove Theorem $1$, we must remove all of these flags from $|\widehat{\mathcal{B}}_2^{II}|$ to obtain $\Delta _1$, using Lemma $2.3$. In a sequence of $35$ steps, we first remove the flags occurring in any orbit after the first orbit listed in the third column, obtaining a subcomplex $\Delta _2$. The remaining $44$ orbits of flags in $\Delta _2$ but not in $\Delta _1$ can easily be removed using Lemma $2.3$, where each flag $\Sigma$ with $x_1=R_p$ is free over a flag $\sigma$, the face of $\Sigma$ with $\sigma$ not containing the group of type $p$.

\medskip
Recall that a maximal simplex $\Sigma$ is {\it free} over some maximal face $\sigma$, if $\Sigma$ is the only maximal simplex with $\sigma$ as a face. In this case we can use Lemma $2.3$ to collapse $\Sigma$ down onto its faces other than $\sigma$. This means that we can remove $\Sigma$ and $\sigma$ from the complex to obtain a homotopy equivalent subcomplex. Note that $\Sigma$ is not necessarily a simplex of maximal dimension; it has to be a maximal simplex, in the sense that it is not a face of a larger simplex in the complex.

\medskip
Since this method depends on the maximality of $\Sigma$ and since the removal of some simplices might make other simplices maximal, the sequence of homotopy retractions we are performing has to be done in the order indicated here.

\medskip
The $35$ steps mentioned above are given in the Tables $3$ and $4$ below. In the first four steps given in Table $3$, we can remove all flags of the stated type. In the next $31$ steps given in Table $4$, we remove only some of the orbits at each stage.

\medskip
We begin by homotopically retracting the pairs of simplices listed in Table $3$.\\

\begin{center}\begin{center}
\begin{tabular}{|l|l|l|}
\hline
Step&$\Sigma$&$\sigma$\\
\hline
$1$.& $p \subseteq M \subseteq {\mathbb L} \subseteq p L \square \subseteq pM L \square$  & $p \subseteq M \subseteq {\mathbb L}  \subseteq pM L \square$\\

$2$.& $M \subseteq {\mathbb L} \subseteq p L \square \subseteq pM L \square$ & $ M \subseteq {\mathbb L}  \subseteq pM L \square$\\

$3$.& $p \subseteq {\mathbb L} \subseteq p L \square \subseteq pM L \square$ & $p \subseteq {\mathbb L}  \subseteq pM L \square$\\

$4$.& ${\mathbb L} \subseteq p L \square \subseteq pM L \square$ & ${\mathbb L}  \subseteq pM L \square$\\
\hline
\end{tabular}
\end{center}

\vspace*{.5cm}
{\bf Table 3:} {\it The first four steps of the homotopy retraction from $|\widehat{\mathcal{B}}_2^{II}|$ to $\Delta_2$}\end{center}

\medskip
We will use the notation: $(k_1) \; x_1 \subseteq \ldots (k_{n-1}) \; x_{n-1} \subseteq x_n$ to denote $k_1 \cdot k_2 \ldots \cdot k_{n-1}$ simplices of the form $x_1 \subseteq \ldots \subseteq x_n$, for which there are $k_i$ choices for the vertices of type $i \in \lbrace 1, \dots n-1 \rbrace$. The notation is needed to record the sizes of the orbits being removed in each step of Table $4$.

\medskip
We will give a detailed description of only one of these steps below. In this description, as well as in the Tables, we abbreviate a flag by using the types for its groups. Since $R_L \not= L$, we abbreviate $R_L$ as the line structure $\mathbb{L}$.

\begin{center}\begin{small}
\begin{center}
\begin{tabular}{|l|l|l|l|}
\hline
Step&$\Sigma$&$\sigma$&Points in $\sigma$\\
\hline

$5$.& $ (8) \; p \subseteq (3)\;  M \subseteq p L \square \subseteq pM L \square$ & $ (24) \; p \subseteq p L \square \subseteq p M L \square$ & $p L \square \setminus p \square$\\

$6$.& $(16)\;p \subseteq p \square \subseteq p L \square \subseteq pM L \square$ & $(16) \; p \subseteq p L \square  \subseteq pM L \square$ & $p \square \setminus \mathbb{L}$\\

$7$.& $(24) \; p \subseteq p \square \subseteq p M \square \subseteq pM L \square$ & $(24) \; p \subseteq pM \square \subseteq pM L \square$ &$ p \square \setminus M_1$\\

$8$.& $(8) \; p \subseteq M \subseteq pM \square \subseteq pM L \square$ & $(8) \; p \subseteq pM \square \subseteq pM L \square$ & $M_1 \setminus p \square$\\

$9$.& $(12) \; p \subseteq (3) \; M \subseteq {\mathbb L} \subseteq p L \square$ & $(36) \; p \subseteq {\mathbb L}  \subseteq p L \square$ &$ \mathbb{L} \setminus L$\\

$10$.& $(8)\; p \subseteq (3)\; M \subseteq p L \square$ & $(24)\; p \subseteq  p L \square$ &$p L \square \setminus p \square$\\

$11$.& $(16) \; p \subseteq  p \square \subseteq p L \square$ & $(16) \; p  \subseteq p L \square$ &$ p \square \setminus \mathbb{L}$\\

$12$.& $(24) \; p \subseteq p \square \subseteq pM \square$ & $(24) \; p \subseteq p M \square$ &$ p \square \setminus M_1$\\

$13$.& $(8) \; p \subseteq M \subseteq  p M \square $ & $(8) \; p \subseteq pM \square$ &$M_1 \setminus p \square$\\

$14$.& $(8) \; p \subseteq  p \square \subseteq p L \square \subseteq pM L \square$ & $(8) \; p \subseteq p \square \subseteq pM L \square$ &$ (\square _2 \cup \square_3) \setminus L$\\

$15$.& $(4) \; p \subseteq (2) \; M \subseteq p L \square \subseteq pM L \square$ & $(8) \; p \subseteq p L\square \subseteq pM L \square$ &$( \square _2 \cup \square_3 )\setminus L$\\

$16$.& $(16) \; p \subseteq p \square \subseteq pM L \square$ & $(16) \; p \subseteq pM L \square$ &$p \square \setminus \mathbb{L}$\\

$17$.& $(8) \; p \subseteq (3) \; M \subseteq pM L \square$ & $(24) \; p \subseteq pM L \square$ &$\mathbb{L} \setminus p \square$\\

$18$.& $(4) \; p \subseteq (2) \; M \subseteq pM L \square$ & $(8) \; p \subseteq  pM L \square$ &$(\square _2 \cup \square _3)  \setminus L$\\

$19$.& $(3) \; p \subseteq (2) \; M \subseteq p L \square \subseteq pM _1 L \square$ & $(3) \; p \subseteq (2) \; M \subseteq pM L \square$ &$L$\\
& \hspace*{1.2cm}$M \in \lbrace M_2, M_3 \rbrace$&&\\

$20$.& $(4) \; p \subseteq p \square \subseteq p L\square \subseteq pM_i L \square$ & $(12) \; p \subseteq p \square \subseteq p L  \square$ &$(\square _1 \cup \square _2 \cup  \square _3) \setminus L$\\
&\hspace*{1.2cm}$i \in \lbrace 1, 2, 3 \rbrace$&&\\

$21$.& $(4) \; p \subseteq M \subseteq p L \square \subseteq pM L \square$ & $(4) \; p \subseteq p L \square \subseteq pM L \square$ & $\square _1 \setminus L$\\

$22$.& $(4) \; p \subseteq p \square \subseteq pM \square \subseteq pM L \square$ & $(4) \; p \subseteq p \square \subseteq pML  \square$ &$\square _1 \setminus L$\\

$23$. & $(4) \; p \subseteq M \subseteq pM \square \subseteq pM L \square$ & $(4) \; p \subseteq pM \square \subseteq pM L \square$ & $\square _1 \setminus L$\\

$24$. & $(2) \;p \subseteq M \subseteq pM \square \subseteq pM L_i \square$ & $(6) \; p \subseteq M \subseteq pM \square$& $\square_1 \setminus p_1$\\
&\hspace*{1.2cm}$i \in \lbrace 1, 2, 3 \rbrace$&&\\

$25$. & $(2) \; p \subseteq p \square \subseteq pM \square \subseteq pM L_i \square$ & $(6)\; p \subseteq p \square \subseteq pM \square$ & $ \square _ 1 \setminus p_1$\\
&\hspace*{1.2cm}$i \in \lbrace 1, 2, 3 \rbrace$&&\\

$26$. & $(2)\; p \subseteq M_1 \subseteq p L \square \subseteq pM_1 L \square$ & $(2) \; p \subseteq M_1 \subseteq pM_1 L \square$ & $ L \setminus p_1$\\

$27$. & $(2)\; p \subseteq p \square \subseteq p L \square \subseteq pM L \square$ & $(2) \; p \subseteq p \square \subseteq pM L \square$ & $L \setminus p_1$\\

$28$. & $(2) \; p \subseteq (3) \; M \subseteq {\mathbb L} \subseteq p L \square$ & $(2) \; p \subseteq (3) \; M \subseteq p L \square$ & $L \setminus p_1$\\

$29$. & $(4) \; p \subseteq (3) \; M \subseteq p L \square$ & $ (12) \; p \subseteq p L \square$ & $ (\square _1 \cup \square _2 \cup \square _3) \setminus L$\\

$30$.&  $(2) \; p \subseteq pM \square \subseteq pM L_i \square$ & $(6) \; p \subseteq pM \square$ & $ \square _1 \setminus p_1$\\
&\hspace*{1.2cm}$i \in \lbrace 1, 2, 3 \rbrace$&&\\

$31$. & $(4)\; p \subseteq M \subseteq pM L \square$ & $ (4) \; p \subseteq pM L \square$&$\square _1 \setminus L$\\

$32$. & $ (2) \; p \subseteq p L \square \subseteq pM L \square $ & $ (2) \; p \subseteq pM L \square$ & $L \setminus p_1$\\

$33$. & $ (2)M \subseteq p L \square \subseteq pM_1 L \square$ & $ (2)\;M \subseteq pM L \square$& \\
&\hspace*{1.2cm}$M \in \lbrace M_2, M_3 \rbrace$&\hspace*{.5cm}$M \in \lbrace M_2, M_3 \rbrace$&\\

$34$. & $(2) \; p \subseteq p \square \subseteq p L_i \square$ & $(30)\; p \subseteq p \square$&$p \square \setminus p_1$\\
&\hspace*{1.2cm}$i \in \lbrace 1, \ldots 15 \rbrace$&&\\

$35$. & $(2)\; p \subseteq {\mathbb L} \subseteq  p L \square$ & $ (2) \; p \subseteq p L \square$ &$L \setminus p_1$\\
\hline
\end{tabular}
\end{center}
\end{small}

\vspace*{.2cm}
{\bf Table 4:} {\it The remaining $31$ steps of the homotopy retraction from $|\widehat{\mathcal{B}}_2^{II}|$ to $\Delta_2$}\end{center}

\newpage
{\it Step 5:} We will discuss this step in detail since the approach we use here will be repeated. Consider a simplex $\Sigma : p \subseteq M \subseteq pL\square \subseteq pML\square$. We consider here the points $p \in pL\square \setminus p\square$. There are $8 \times 3$ such points, $8$ points in each of the three $\mathcal{M}$-spaces of $pML\square$. This means that for a fixed $pL\square \subseteq pML\square$, there are $24$ simplices $\Sigma$ of the type chosen above. We write this as $(8)\; p \subseteq (3) \; M \subseteq pL\square \subseteq pML\square$. Fix one such simplex $\Sigma$. This is a simplex of maximal dimension and free over the face $\sigma : p \subseteq pL \square \subseteq pML\square$. To see this note that $ p \not \in L$. Thus $\langle p, L \rangle$ is a plane on $L$ and thus lies in a unique $\mathcal{M}$-space. Note the simplex $p \subseteq \mathbb{L} \subseteq pL\square \subseteq pML\square$ was removed in Step $3$. Thus $\Sigma$ is free over $\sigma$. Finally we can apply Lemma $2.3$ and remove these two simplices.

\medskip
The remaining steps are listed in Table $4$. In the first column we number the steps. In the second column we give the maximal simplex $\Sigma$ and in the third column we give a face $\sigma$ of $\Sigma$ such that $\Sigma$ is free over $\sigma$. In the fourth column we specify the collection of points in $\sigma$ being removed. We use the notation described above.

\medskip
After these steps we are left with the complex $\Delta _2$. There are still $44$ orbits of flags remaining in $\Delta_2 \setminus \Delta_1$, all of these corresponding to the first orbit listed in the third column of Table $2$. But these can easily be canceled in pairs, with a flag $\Sigma$ having $x_1=R_p$ of type $p$ being free over its face $\sigma$ which does not involve the group of type $p$. Recall that by Proposition $5.1, \; \Delta _1$ is homotopy equivalent to $\Delta_0$ which, in its turn, is homeomeorphic to $\Delta$. This ends the proof of Theorem $1$.\\

\section{On the Lefschetz module}

Let $G = Co_3$ and let $\Delta$ be the $2$-local geometry of $G$ described in Section $3$.

\begin{prop}Let $z$ be a central involution in $G$. The set $\Delta^z$ is contractible.
\end{prop}

\begin{proof}Recall that $C_G(z) = G_p$ for the point $p=\langle z \rangle$. It follows that $p *
\text{Res}(p)=\text{Star}(p) \subseteq \Delta^z$. It is easy to see that a point $q \in \Delta ^z$ if and only if $q \in p ^{\perp}$. The main idea of the proof is to use Lemma $2.2$ to homotopically remove those lines and $\mathcal{M}$-spaces from $\Delta^z$ that do not contain $p$, as well as those points in $p ^{\perp}$.

\medskip
It follows from $(\Delta 4)$, Section $3$, that if a line $L$ is an element of $\Delta^z$ then $p$ is collinear with one or all the points of
$L$. Thus there are three types of lines we have to consider: lines on $p$, already included in $\text{Res}(p)$; lines with only one point fixed, the other two points being
interchanged; lines with all three points fixed. By $(\Delta 5)$, given a point, $\mathcal{M}$-space pair $(p, M)$, the set $p^{\perp} \cap
M$ is at most one line or $p \in M$. This implies that an
$\mathcal{M}$-space $M \in \Delta ^z$ can have one point in $\Delta^z$, three points forming a line in $\Delta^z$ or $p \in M$. This is easily seen if we
recall that $M$ contains $15$ points, thus under the action of $z$, at
least one point is fixed.

\medskip
Let us consider the case when the line $L$ does not contain $p$ and has
only one point $q$ fixed by the action of $z$. Then the collection
$\text{Res}_{\Delta^z}(L) = \lbrace q \rbrace * \lbrace \text{some}\; \mathcal{M}-
\text{spaces on}\; L \rbrace$ is a cone on $q$, and thus contractible. It follows that we can remove
all such lines from $\Delta ^z$, yielding a homotopy equivalent
subcomplex $\Delta ^z_1 \subseteq \Delta ^z$.

\medskip
Let $M \in \Delta_1 ^z$ be such that $p^{\perp} \cap M=L$, a line: this means all
three points of $L$ are in $\Delta_1^z$. Then $\text{Res}_{\Delta_1^z}(M)$
consists of $L$ and the three points on $L$, which is a cone on $L$ and thus contractible. The other lines of $M$ are either not stabilized by $z$ or were removed at the previous step, so are not in $\Delta_1^z$. We can remove all such
$\mathcal{M}$-spaces from $\Delta_1 ^z$, yielding a homotopy equivalent
complex $\Delta ^z_2 \subseteq \Delta ^z_1 \subseteq \Delta ^z$.

\medskip
Let us now consider a line $L$ whose three points are all fixed by the action of $z$. Then $p$ and $L$ generate a plane $\langle p, L \rangle$
which lies in a unique $\mathcal{M}$-space $M$. The other
$\mathcal{M}$-spaces from $\Delta ^z$ through $L$ were deleted at the
previous step. Thus $\text{Res}_{\Delta ^z_2}L=\lbrace M \rbrace * \lbrace 3\;
\text{points on} \;L \rbrace$, again a cone on $M$ and thus contractible. Remove
these lines, yielding $\Delta ^z_3 \subseteq \Delta ^z_2$. Now all the lines of $\Delta ^z$ which are not on $p$ have been removed.

\medskip
Let $M \in \mathcal{M}$ be such that $p^{\perp} \cap M=\lbrace q
\rbrace$, a single point. Then $\text{Res}_{\Delta ^z_3}M$ is the point $q$
itself. Remove these $\mathcal{M}$-spaces, yielding $\Delta ^z_4
\subseteq \Delta ^z_3$.

\medskip
Next, consider the collection of points collinear with $p$ and look at
their residues in $\Delta ^z_4$. Let $q \in p^{\perp}$. Since we have
already contracted all the lines on $q$ which do not contain $p$, the residue $\text{Res}_{\Delta_4 ^z}(q)$ contains a unique line $\langle p, q \rangle$ and the three $\mathcal{M}$-spaces on that line, therefore this residue is contractible. So we can remove these points to retract $\Delta_4 ^z$ to $\text{Star}(p)$. Of course, $\text{Star}(p)$ is contractible, a cone on $p$.
\end{proof}

\bigskip
In what follows we will use the following results due to Th\'{e}venaz; see \cite{th87}, Theorem $2.1$ and Corollary $2.3$. A group $B$ is called cyclic mod $p$ if the quotient group $B/O_p(B)$ is cyclic.

\begin{thm}\cite[Th\'{e}venaz]{th87} a). Let $\mathfrak{X}$ be a class of subgroups of $A$ which is closed under subconjugacy. Let $\Delta$ be an admissible $A$-complex such that the reduced Euler characteristic $\tilde \chi (\Delta ^B) = 0$ for every subgroup $B$ which is cyclic mod $p$ and satisfies $O_p(B) \not \in \mathfrak{X}$. Then the reduced Lefschetz module $\tilde L (\Delta)$ is a $\mathbb{F}_p A$-virtual module projective relative to the collection $\mathfrak{X}$.\\
b). Let $p^n = |A|_p$ be the $p$-part of $A$,  and let $p^k$ be the highest power of $p$ dividing the order of some subgroup belonging to $\mathfrak{X}$. Then $\tilde \chi (\Delta)$ is a multiple of $p^{n-k}$.
\end{thm}

We let $\mathfrak{X} = \mathcal{B} _2 \setminus \mathcal{\hat B}_2$, the family of the radical $2$-subgroups of $G=Co_3$ which are pure non-central; these are subgroups in the conjugacy classes of $\lbrace R_2, R_3, R_4 \rbrace$. It is clear that $\mathfrak{X}$ is closed under subconjugacy.

\begin{rem}The subgroups $Q$ in $G=Co_3$ which are cyclic mod $2$ and such that $Z(O_2(Q))$ contains no central involutions have $O_2(Q)$ in $\mathfrak{X}$. To see this, let $Q$ be a cyclic mod $2$ subgroup of $G$; set $R = O_2(Q)$ and $H = \Omega _1 Z(R)$. Thus $R$ is contained in a member of $\text{Syl}_2(N_G(H))$. Let us assume $H$ contains no central involutions; then $H$ is of the form $2, 2^2$ or $2^3$; see \cite[Section 5]{fin73}. The normalizers of $R_2, R_3$ and $R_4$ are given in Table $1$, at the beginning of Section $4$. It is easy to see from Table $1$, that the Sylow $2$-subgroup of $N_G(H)$ has the form $H \times K$ with $K$ some $2$-subgroup. If $H \leq R \leq H \times K$ then $R = H \times K_1$ for some subgroup $K_1 \leq K$. But $Z(R) = H \times Z(K_1)$, where $Z(K_1)$ is nontrivial if $K_1$ is nontrivial. Thus $\Omega _1 Z(R) = H \times \Omega _1 Z(K_1) = H$. It follows that $\Omega _1 Z(K_1) = 1$ and also $Z(K_1) = 1$ which, in its turn implies $K_1 = 1$. Finally $R = H$.
\end{rem}

Furthermore, $G$ acts admissibly on $\Delta$ and $\tilde \chi (\Delta ^z) = 0$, for every central involution $z \in G$, according to Proposition $6.1$. Now, by P.A. Smith theory, it follows that $\tilde \chi (\Delta ^B) = 0$ also, for every subgroup $B$ which is cyclic mod $2$ and such that $O_2(B) \not \in \mathfrak{X}$. Thus, all the necessary conditions required by the Theorem $6.2$ are in place and we can formulate the following:

\medskip
\begin{tb} {\it Let $G = Co_3$ and let $\mathfrak{X} = \mathcal{B}_2 \setminus \mathcal{\hat B}_2$. Let $\Delta$ denote the $2$-local geometry of $G$. The reduced Lefschetz module $\tilde L (\Delta)$ is a virtual $\mathbb{F}_2G$-module projective relative to the collection $\mathfrak{X}$. Furthermore $\tilde \chi (\Delta)$ is a multiple of $2^7$.}
\end{tb}

\medskip
\begin{rem} Note that $|G|_2 = 2^{10}$ and that the $2$-part of $\tilde \chi (\Delta)$ is $2^7$. Also, the largest subgroup in $\mathfrak{X}$ has order $2^3$. Th\'{e}venaz' result tells us that the reduced Euler characteristic is a multiple of $2^7$, without giving an upper bound.\\
\end{rem}

{\bf Acknowledgements.}

\medskip
We would like to express our gratitude to S.D. Smith for suggesting the problem and for helpful comments. We would also like to thank the referee whose suggestions led to an improved exposition of the paper.\\

\end{document}